\documentclass[a4paper,12pt]{article}
\usepackage{latexsym}
\usepackage{amsfonts}
\usepackage{amsmath}
\def \E {{\rm E }}
\def \R {{\rm I\!R}}

\def \N {{\rm I\!N}}
\def \t {{\tau}}

\def \grad {{\rm grad}}
\def \< {\left< }
\def \> { \right> }

\title{Propagation of microlocal singularities for stochastic partial differential equations} 

\author{ Adnan Aboulalaa}
\date{}

\begin{document}
\vskip2cm
\newtheorem{thm}{Theorem}[section]
\newtheorem{defi}[thm]{Definition}
\newtheorem{prop}[thm]{Proposition}
\newtheorem{Rq}[thm]{Remark}
\newtheorem{lm}{Lemma}[section]
\newtheorem{co}[thm]{Corollary}
\maketitle

\footnotetext[1]{ 
Email: adnan.aboulalaa@polytechnique.org.}
 
\begin{abstract}
Microlocal analysis techniques are extended and applied to stochastic partial differential equations (SPDEs). 
In particular, the H\"ormander propagation of singularities theorem is shown to be valid for hyperbolic SPDEs driven by a standard Brownian motion. In this case the wave front set of the solution is invariant under the stochastic Hamiltonian flow associated to the principal symbol of the SPDE. This study leads to the introduction of a class of random pseudodifferential operators.
\end{abstract}

\section{Introduction}

The singularities of the solutions to linear partial differential equations (PDEs) have been studied for a long time and these studies culminated 
in the early 1970s with microlocal analysis; this framework has provided a refined mathematical formulation of the concept of singularities and the corresponding propagation theorem is one of the remarkable results obtained in this context.

\noindent
When the partial differential equations govern some systems involving random media or random perturbations, one is led to consider stochastic partial differential equations (SPDEs) models, a subject that has been extensively studied during the last three decades.  

\noindent
The investigation of singularities of SPDEs together with the possible applicability of microlocal analysis techniques is a natural issue in this perspective. Several problems may be addressed: for what type of SPDEs may these techniques be relevant ? how to deal with randomness when considering pseudodifferential operators ? is it possible define singularities in a way similar to the deterministic linear PDEs and is there a propagation phenomenon and theorem for SPDEs ? etc.

The purpose of this paper is to address such issues for hyperbolic stochastic partial differential equations of the type:
\begin{equation}
\nonumber
(E)\left\{ \begin{array}{l}
  \displaystyle
 du(t)=\sum_{i=1}^{n}[a_{i}(t,x,D)u(t)\circ dw^{i}(t)+f_{i}(t)\circ dw_{i}(t)]+b(t,x,D)u(t)dt+g(t)dt,\\
\displaystyle
 u(0)=u_{0}\in (H^{s}(\R^{d}))^{d'},
\end{array}\right.
\end{equation}
This is the stochastic counterpart of linear hyperbolic systems which were studied in particular by Friedrichs \cite{Friedrichs}. These systems give
rise to a large class of hyperbolic PDEs, such as second order linear hyperbolic PDEs, Maxwell equations, Dirac equations, etc.
The existence and uniqueness of the solution to the stochastic hyperbolic systems $(E)$ were studied in \cite{Ab1}, where $a_{i}(t,x,D),b(t,x,D)$ are smooth families of $d'\times d'$-matrices of first order pseudodifferential operators (PDO), $w^{i}(t), t\in I$ are standard Wiener processes, $f_{i},g$ are continuous, possibly random, functions from $I=[0,T], T>0$ to $(H^{s}(\R^{d}))^{d'})$, $\circ$ corresponds to the Fisk-Stratonovich integral or differential, and $H^{s}(\R^{d})$ is a Sobolev space $(s\in \R, d, d'\geq 1 )$.

There has been a wide interest in the topic of singularities of PDEs solutions during the 20th century in accordance with the concepts and techniques available. 
At the beginning, singularities were associated with the discontinuities of the solutions or their derivatives and were studied for wave and Maxwell equations. 
It was shown that the discontinuities of the solutions of wave equations in inhomogeneous media propagate along the bi-characteristic curves which correspond to the paths given by 
the geometrical optics, see, e.g., Luneberg \cite{Luneburg1}, Kline\cite{Kline}.
This fact was important not only because it establishes a formal link between wave and geometric optics, viewing the later as an approximation of the former and
providing a tool to justify the practical calculations based on geometrical optics, but also because of its conceptual side: this correspondence can be applied to other areas
such as acoustics (see \cite{Keller1}) and mechanics.
This result was subsequently generalized to other types of PDEs, possibly non-linear, see Courant and Hilbert \cite{CH}, and John \cite{John}. It was also extended to
linear symmetric systems (Courant and Lax \cite{CL}).
\\
On the other hand, the fact that geometrical optics (Fermat's principle, eikonal equation) forms an approximation of wave optics was noticed before from a  
different point of view, that is to show that the solutions to geometrical optics equations correspond to those of wave optics when the wave lengths become very small 
(Runge - Sommerfeld \cite{RS}). Let us recall in passing that this fact was actually a 
fundamental starting point of quantum mechanics in Shr\"odinger's work: classical mechanics whose equations are similar to those of geometrical optics, may be considered as
an approximation of a wave mechanics for which we have to find the equation, and based on de Broglie's work, Shr\"odinger was able to find his famous equation. 
This point a view was subsequently clarified, starting with Birkhoff's work \cite{Birkhoff1} and extensively studied in the framework of semi-classical approximation or analysis (see, e.g., 
\cite{Martinez}, \cite{Zworski}, \cite{GS})
\\
The description of singularities was refined in the late 1960s by Sato and H\"ormander with the concept of singular spectrum or wave front set of a distribution; this set $WF(u)$ provides the couples $(x, \xi)$ such that $u$ is 
singular in the point $x$ and in the frequency direction given by $\xi$, which means that $\hat{u}$ is not rapidly decreasing in the direction $\xi$. In fact the thus defined
wave front provides the space points $x$ where $u$ is not smooth, which are viewed as singularity points whereas the frequency direction $\xi$ can give information about
the form of the set of singularities or the direction of propagation (see $\S 2 $ below), whence a link with the intuitive term wave front used by the physicists to refer to wave propagation phenomena. With this concept, 
H\"ormander obtained the propagation of singularities theorem of the solution to the hyperbolic equation $\partial_{t} u= A(x,D)u$, where $A(x,D)$ is a first order pseudodifferential operator with a classical symbol; it states that the wave front set of the solution $u(t, \cdot)$ is invariant under the bicharacteristic curves $\chi_{t}$ of the principal symbol $a_{1}$ of the operator $A$, i.e. the integral curves of the flow $(\grad_{\xi} a_{1}, - \grad_{x} a_{1})$, see H\"ormander \cite{H4}. This result was also established for more general equations of the form $P u = 0$.
\\
The goal of this paper is to show that this result is still valid for the hyperbolic stochastic differential equations $(E)$, that is $WF(u(t, \cdot))=\chi_{t} WF(u(0, \cdot))$,
where this time $\chi_{t}$ is the stochastic integral Hamiltonian flow associated to the principal symbols of the operators $a_{i}(x,D)$ which give the stochastic part with
 Stratonovich differential equations and $b(x,D)$ which gives the deterministic part.
\\
For this purpose, we shall be led to introduce a class of random pseudodifferential operators which is more general than the ones used in some previous works 
(see Dedik, Shubin \cite{DS2}, Pankov \cite{Pankov}, Fedosov-Shubin \cite{FS}). Some properties of these random PDOs will be established; let 
us mention that Liu and Zhang \cite{LZ} studied a class of random PDO close to the one considered in this paper.

In the case of SPDEs, hyperbolic equations are often studied with a space-time white noise, in which case the solutions are at most H\"older continuous. Despite this lack
of regularity, the problem of propagation of singularities was considered by Walsh (\cite{Walsh1}, \cite{Walsh2}) for the Brownian sheet which
actually corresponds to the solution of the wave equation with a space-time white noise source; this study was extended to a nonlinear stochastic wave equation 
by Carmona and Nualart\cite{CN}. In this case singularities mean a failure to have some modulus of continuity. See also Blath and Martin \cite{BM} for an extension of this kind of study to semi-fractional Brownian sheets.
\\
This paper is organized as follows: Notations and preliminaries on linear hyperbolic SPDEs are presented in section 2. In order to position the problem addressed in 
this paper, it is relevant to revisit the results so far obtained on propagation of singularities for PDEs and SPDEs; section 3 provides a quick overview 
of these results. The notion of random pseudodifferential operators is introduced in section 4, with some results on the stochastic integration of random symbols that will be needed later. The definition of these random PDOs requires a pathwise bound of the derivatives of random symbols. 
One of the complications encountered is the verification of these pathwise bounds when the random symbols used are defined by stochastic integrals. 
This difficulty is overcome with the use of a device based on the Kolmogorov-Centsov continuity criterion. 
Section 5 contains the statement and proof the main result, that is the propagation of singularities of the hyperbolic SPDEs $(E)$ along its stochastic Hamilton flow. The proof is based on the pathwise approximation of the solutions established in \cite{Ab1} and the properties of random symbols obtained in section 4.


\section{Notations and preliminaries}
\subsection{Notations and settings}
We consider SPDEs of type (E), driven by a finite dimensional Brownian motion. Let $(w_{t}), t\in I:=[0,T]$ be a one-dimensional
 Brownian motion  defined on a filtered probability space $(\Omega,{\cal F}, {\cal F}_{t}, P)$ with ${\cal F}_{t}=\sigma(w_{\t}, \t\leq t)$ and $T>0$ fixed throughout this paper. 
We also fix $d,d'\geq 1$. The solutions to (E) will be considered as $(H^{s})^{d'}$-valued process, where $H^{s}$ is the Sobolev space of order $s$;
$\left< \cdots, \cdots \right>_{s}$ will be the scalar product on the Sobolev space $H^{s}:=H^{s}(\R^{d}), s\in \R$ and the same notation is used for the scalar
product on $(H^{s})^{d'}$. We refer to \cite{Ab1} for a detailed study of the stochastic hyperbolic systems (E) and their properties.  

For $p>0$, $M^{s}_{p}(I, (H^{s})^{d'})$ denotes the  set of adapted $(H^{s})^{d'}$-valued processes $u(t),t\in I$ such that 
$\|u\|_{s,p}:=(E\sup_{t\leq T}|u(t)|_{s}^{p})^{1/p}< +\infty$, this quantity being the norm of $u$ in $M^{s}_{p}$ 
\\
We use the standard PDO notations: for $\alpha = (\alpha_{1}, \cdots, \alpha_{d}) \in \N^{d}$, $|\alpha|= \sum_{j=1}^{d}\alpha_{j}$ and $D^{\alpha}= D^{\alpha_{1}}_{1} \cdots D^{\alpha_{d}}_{d}$, where
$D^{\alpha_{j}}_{j} = (-i)^{\alpha_{j}}\partial_{j}^{\alpha_{j}}$ with $i=\sqrt{-1}$ and $\partial_{x_{j}}=\partial_{x_{j}}$, for $x=(x_{1}, \cdots, x_{d})$. 

As usual $S^{m}$ is the set of symbols $a(x,\xi)$ of order 
$m$ on $\R^{d}$, i.e. $a\in C^{\infty}(\R^{d}\times \R^{d})$ and
for all $\alpha, \beta \in \N^{d}$ there is a constant $C(\alpha,\beta)$ 
such that $|D_{\xi}^{\alpha}D_{x}^{\beta}a(x,\xi)|\leq C(\alpha, \beta)
(1+|\xi|)^{m-|\alpha|}$. 
For such a symbol, $a(x,D)$ denotes the associated pseudodifferential operator defined by $a(x,D)u(x)=\int a(x,\xi) \hat{u}(\xi)
e^{i x.\xi}d\xi $ for $u\in C^{\infty}_{0}(\R^{d})$, $a^{*}(x,D)$ is the adjoint of
 $a(x,D)$ and ${\rm OPS^{m}}$ will denote the set of such operators.
\noindent
As we consider symmetric systems, the solutions are vector valued (in $\R^{d'}$) and the operators $a(x,D)$ in (E) are matrices of pseudodifferential 
operators: $a(x,D)=(a^{ij}(x,D), i,j=1,...,d')$
with $a^{ij}(x,D)\in {\rm OPS}^{m}$ for some $m$; the set of such matrices of PDOs is also denoted by ${\rm OPS}^{m}$.
\\
\\
Let $X\subset\R^{d}$ be an open set and $u\in {\cal D}'(X)$. We denote by $WF(u)\subset X \times \R^{d}-\{0\}$
 the wave front set of the distribution $u$ defined by: $(x_{0},\xi_{0})
\notin WF(u)$ if there is a $\phi\in C_{0}^{\infty}(\R^{d})$ with $\phi\equiv
 1$ near  $x_{0}$ such that $\widehat{\phi u}(\xi)$ is rapidly decreasing as
$|\xi|\rightarrow \infty$ in an open cone $\Gamma$ containing $\xi_{0}$; $\Gamma \subset X \times \R^{d}$ is said
to be a cone if for all $(x, \xi)\in \Gamma$ we have $(x, t\xi)\in \Gamma, \forall t\in ]0, +\infty]$.
\\
The wave front (\cite{H1}), also called the singular spectrum of a distribution $u\in {\cal D}'(X)$ is a refinement of the singular 
support ${\rm sing\; supp} (u)$ defined by $x_{0} \notin {\rm sing\; supp} (u)$ iff there exists a $\phi\in C_{0}^{\infty}(\R^{d})$ with $\phi (x_{0})\neq 0$ and
$\phi u \in C_{0}^{\infty}(X)$. We have indeed: ${\rm sing \; supp} (u) = {\rm Pr}(WF(u))$ where ${\rm Pr}: (x, \xi) \mapsto x$ is the projection on the $x$-space $X$.
\\
The wave front set provides the frequency directions $\xi$ in which the singularities occur. At a first sight, there is no direct link with the physical intuitive notion of wave fronts which gives the form of curve or surface of singularities and the direction of their propagation. 
However these elements of the frequency space, the $\xi$, give indeed such information: for instance let $u(x_{1}, x_{2})$ be a function on $\R^{2}$ which is smooth except a singularity at the line $x_{1}=a$; then its wave front is the set $\{ ((a,x_{2}), (\xi_{1}, 0)): x_{2} \in \R, \xi_{1} \in \R^{*} \}$: one can indeed see that:

$$\widehat{\phi u}(0,\xi_{2})= \int e^{i x_{2} \xi_{2}}d x_{2} ( \int u(x_{1},x_{2}) d x_{1}),$$

\noindent
and from the assumptions on $u$ it follows that $ \int u(x_{1},x_{2}) d x_{1}$ is smooth as 
a function of $x_{2}$, so that $\hat{\phi u}(0,\xi_{2})$ is rapidly decreasing as a function of $\xi_{2}$. In this example the direction $(\xi_{1}, 0)$ corresponds indeed to 
vectors perpendicular to the set of singularities $x_{1}=a$ and hence provides an information about the direction of singularities. This elementary example can be extended to more general situations with singularities given by curves or surfaces and with a time variable representing the evolution which will be used to see the direction of propagation of singularities.
\\
Equivalently, we have the following characterization of the wave front set, which is usually used in the study of singularities and their propagation:
$WF(u)=\bigcap\{ {\rm Char} (p): p(x,D)u \in C^{\infty}, p\in S^{0}_{ph} \}$, where 
$S^{0}_{ph}$ is the set of polyhomogeneous symbols in $S^{0}$, that is $p(x,\xi)\sim
 \sum_{i=0}^{\infty}p_{j}(x,\xi)$ with $p_{j}(x,\xi)$ homogeneous of degree
$-j$ in $\xi$ for $|\xi|>1$ and the charateristic set ${\rm Char}(p)$ is defined by ${\rm Char}(p)=\{(x,\xi): p_{0}(x,\xi)=0 \}$. 
\\
\noindent
This means that $(x,\xi)\notin WF(u)$ iif we can find $p\in S^{0}_{ph}$ with $p_{0}(x,\xi)\neq
0$, and $p$ vanishes outside some conic neighborhood of $(x,\xi)$ and 
$p(x,D) u \in C_{0}^{\infty}$. 


\section{Propagation of singularities for PDEs and SPDEs}

\subsection{The case of partial differential equations}

Before the development of microlocal analysis, the singularities of a solution $u$ of a PDE meant discontinuities of $u$ or its first or higher order derivatives. The corresponding results on propagation of singularities concerned thus the time evolution of these  discontinuities. We start this overview by a digression on the arguments used to show how the bicharacteristic curves carry the discontinuities for hyperbolic equations (see Luneburg \cite{Luneburg1}, Courant, Hilbert \cite{CH}, John \cite{John}, Lax \cite{Lax2}). We consider the simplest hyperbolic equation in a two dimension space:
\begin{equation}
\label{review1}
 a(x,y) \partial_{x}^{2} u + b(x,y) \partial_{xy}^{2} u+ c(x,y) \partial_{y}^{2} u +d(x,y)=0
\end{equation}
Let ${\cal C}$ be a curve defined by $\phi(x,y)=0$ which divides the plane $x,y$ into two regions $ O_{-},  O_{+}$ where $\phi < 0 $ and $\phi > 0 $ respectively. We
suppose that the derivatives $\partial_{x} u, \partial_{x} u $ are continuous on ${\cal C}$ but the 2d derivatives have jumps
across this curve: we suppose that $u$ is $C^{2}$ in the region $ O_{-},  O_{+}$. If $s$ is a parametrization of the curve 
${\cal C}_{c}$: $\phi(x,y)= c= {\rm Constant}$, that is: $\phi(x(s),y(s))=c$ then we have $\phi_{x}(x(s), y(s)) x'(s) + \phi_{y}(x(s), y(s)) y'(s)=0 $ and by a suitable choice of $s$ the interior derivatives along a curve ${\cal C}_{c}$, which we denote by $u_{x}'$, is given by: $u_{x}'= u_{xx} \phi_{x}- u_{xy}\phi_{y} $; this expression is taken at the points $(x(s), y(s))$. Similarly we have $u_{y}'= u_{yx} \phi_{x}- u_{yy}\phi_{y} $. If we suppose that $u_{x}', u_{y}'$ are continuous, then their jumps $[u_{x}'], [u_{y}']$ across ${\cal C}$ will be $=0$ and we will have:
\begin{equation}
\label{review2}
   [\partial_{x}^{2}u] \partial_{y} \phi - [\partial_{xy}^{2}u] \partial_{x} \phi=   [\partial_{xy}^{2}u] \partial_{y} \phi - [\partial_{yy}^{2}u] \partial_{x} \phi=0
\end{equation}
and $ [\partial_{yy}^{2}u] \partial_{x} \phi - [\partial_{yx}^{2}u] \partial_{y} \phi$, which implies that:
\begin{equation}
\label{review3}
   [\partial_{x}^{2}u] = k (\partial_{x} \phi)^{2}, \;   [\partial_{y}^{2}u] = k (\partial_{y} \phi)^{2}, \;   [\partial_{xy}^{2}u] = k \partial_{x} \phi \partial_{y} \phi
\end{equation}
Now let $P\in {\cal C}$ and $P_{1} \in O_{-}, P_{2} \in O_{+}$ be 3 points close to each other and consider the value of a solution to (\ref{review1}) 
when $P_{1}, P_{2}$ tend to $P$; then we will have:
\begin{equation}
\label{review4}
   a [\partial_{x}^{2}u] + b [\partial_{xy}^{2}u] +c [\partial_{y}^{2}u] =0
\end{equation}
Using (\ref{review3}) we see that $\phi$ must satisfy the following equation:
\begin{equation}
\label{review5}
   a \partial_{x}\phi^{2} + b \partial_{x}\phi \partial_{y}\phi +c \partial_{y}\phi^{2} =0
\end{equation}
which is the characteristic equation of (\ref{review1}). From this we deduce that if we call ``wave front surface'' the surface ${\cal C}$ across which the solution $u$ to
 (\ref{review1}) has discontinuities for some higher order derivatives (and ${\cal C}$ separates the space into two regions on which the solution $u$ is regular, see Luneberg \cite{Luneburg1}, \cite{CH}), then this wave front must be a characteristic surface of (\ref{review1}). Now since (\ref{review5}) is a 1st order partial differential equation, 
one can see that the discontinuities are carried on the integral curves of this PDE, which are called characteristics, and since the PDE (\ref{review5}) is already a characteristic equation, these curves are called the bicharacteristics of the PDE (\ref{review1}), see \cite{Luneburg1}. In the case of a wave equation where $y=t, a=1/c^{2}, b=0, c=1$, we recover the eikonal equation:
\begin{equation}
\label{review6}
   \partial_{t}\phi^{2}(t,x) -\frac{1}{c^{2}} \partial_{x}\phi^{2}(t,x) =0
\end{equation}
The property of propagation of discontinuities along the bicharacteristics was extended by Courant and Lax to the hyperbolic systems $\sum_{j=1}^{p} A_{j}(x) \partial_{j} u + B(x) u= 0 $ with some additional conditions.
Luneburg \cite{Luneburg1} was probably the first to study this propagation of discontinuities in a clear framework for the Maxwell equations, considering the rays as the orthogonal curves to characteristic surfaces. He showed that they correspond to the geometrical optics rays which are given by Fermat's principle. 
This property represents the second aspect of the links between wave and geometrical optics.
\\
The first aspect of this connection is the fact that the solutions of wave optics equations tend to those of geometrical optics when the wave
lengths become small or, equivalently, in the high frequency limit. This remark went back to the works of Kirchhoff \cite{Kirchhoff}, Debye, Sommerfeld and Runge \cite{RS}. 
We recall that the idea of this approach is to express the solution to an equation of the type:
\begin{equation}
\label{review61}
   \frac{\partial^{2} u }{\partial t^{2}}(t,x) -\frac{1}{c^{2}(x)} \frac{\partial^{2} u }{\partial x^{2}}(t,x) =0
\end{equation}
as: $u(t,x)=e^{i\omega t} A(x)$; here the velocity can be written as $c(x)=c_{0}/n(x)$, where $c_{0}$ is the velocity in the vacuum or in an homogeneous medium and $n(x)$ is
the index of the inhomogeneous medium. Then $A(x)$ verifies:
\begin{equation}
\label{review7}
  A''(x) + k^{2}(x) A(x) =0 \; \; {\rm with} \; k(x)= \frac{\omega n(x)}{c_{0}}
\end{equation}
and if we consider a high frequency expansion of the form:
\begin{equation}
\label{review8}
  A(x)=  e^{i\omega \phi(x)} (A_{0}(x)+\frac{1}{\omega}A_{1}(x)+\frac{1}{\omega^{2}}A_{2}(x)+ \cdots ),
\end{equation}
then the phase $\phi(x)$ will satistfy the eikonal equation:
\begin{equation}
\label{review81}
 |\frac{d \phi}{dx}(x)|= \frac{n(x)}{c_{0}},
\end{equation}
 while the amplitudes verify transport equations of the type: $\phi''(x) A_{0}(x)+2\phi'(x) A_{0}'(x)=0$, etc. This point of view was also studied in order to give a
satisfactory mathematical basis to the assertion that geometrical optics corresponds to a high frequency limit of wave optics. The same idea was applied to mechanics in order
to establish that the limit as the Plank constant $h\longrightarrow 0$ of the quantum mechanical description (via Shr\"odinger equation) corresponds to the classical mechanic equation (Birkhoff \cite{Birkhoff1}). This high frequency analysis was extensively studied afterward within the framework
of semi-classical approximations.
The two points of views (propagation of singularities and high frequency limit) have received modern formulations through microlocal analysis, see Garding \cite{Garding} for a detailed account of this subject. Hormander's theorem on propagation of singularities appeared in several papers (see, e.g., \cite{H2}, \cite{H3}, \cite{DH}); for the case of
hyperbolic PDEs it is formulated as follows: let $u$ be the solution to:
\begin{equation}
\label{review82}
 \frac{\partial u }{\partial t}= a(t,x, D) u + f, \; \; u(0,x)=u_{0} \in H^{s}
\end{equation}
where $a(t,x, D)$ is a first order family of PDOs whose symbols $a(t,x,\xi)$ form a bounded family in $S^{1}$, with $t\mapsto a_{t}$ continuous and 
$ a(t,x,\xi)\sim \sum_{j=0}^{\infty}a_{1-j}(t,x,\xi) \in S^{1}_{ph}$. Then we have $WF(u(t, \cdot))=\chi_{t} WF(u(0, \cdot))$, that is, the wave front of $u(t,0)$ is invariant
under the Hamiltonian flow $\chi_{t}$ of the principal symbol of $a(t,x,D)$ given by:  $\chi_{t}(x, \xi)=(x_{t}, \xi_{t})$ with:
\[  
\frac{d x_{t}}{dt}=-\frac{1}{i}\frac{\partial a_{1}}{\partial\xi_{i}}(t,x(t),\xi(t))\frac{\partial}{\partial x^{i}},  
           \; \; \;    \frac{d \xi_{t}}{dt} = \frac{1}{i}\frac{\partial a_{1}}{\partial x^{i}}(t,x(t),\xi(t))\frac{\partial}{\partial \xi_{i}}
\]
and $x_{0}=x, \xi_{0}=\xi$. This theorem was extended to other PDEs of the form $P u =f$ where $P$ is a pseudodifferential operator of order $m$ such that its 
principal symbol $p_{m}$ is homogeneous of degree $m$; it is also assumed that $P$ is properly supported and $p_{m}$ is real. 
Let 
\[H_{p_{m}}=\sum_{j=1}^{d} \frac{\partial p_{m}} {\partial {x_{j}}} \frac{\partial } {\partial {\xi_{j}}} 
                          - \frac{\partial p_{m}} {\partial {\xi_{j}}} \frac{\partial } {\partial {x_{j}}}
													\]
be the Hamiltonian vector field of $p_{m}$; its integral curves $\gamma_{t} (x_{0}, \xi_{0}) =(x_{t}, \xi_{t})$ are called the bicharacteristics 
and among these curves, those for which
we have: $p_{m}(x_{t}, \xi_{t})=0$ are called the zero bicharacteristics. The propagation of singularities theorem in this case asserts that if  $P u =f$ with $f\in {\cal D}'(U)$, $U\subset \R^{d}$ being an open set, then: $WF(u) - WF(f)\subset {\rm Char} (P)$ and is invariant under the zero bicharateristics. In particular, if $f$ is smooth then
$ (x_{0}, \xi_{0}) \in WF(u)$ iif $ (x_{t}, \xi_{t}):= \gamma_{t}(x_{0}, \xi_{0}) \in WF(u)$. This result has been reformulated and studied in several situations, including
the problems of reflection, diffraction, etc., see, e.g., \cite{Taylor1}, \cite{Raltson}.

\subsection{The case of stochastic partial differential equations}

In the case of SPDEs, we have two situations: in the first one, these equations are driven by standard or cylindrical Brownian motion and the solutions may have sufficient regularity properties, as being $C^{\infty}$ or being elements of Sobolev spaces $H^{s}$. 
In this case one can hope for a propagation of singularities result similar to
deterministic PDEs, provided we deal with hyperbolic type equations, since the parabolic type equations do not lead to propagation phenomena. This is precisely the purpose of the results presented in the next sections.
\\
In the second situation, one considers SPDEs driven by a space-time white noise, in which case the solutions are not regular and they are at most H\"older continuous.  Singularities are then defined as a failure to have a local modulus of continuity, which is related to the law of iterated logarithm (LIL) for the models so far studied.
More precisely, Walsh \cite{Walsh1}, \cite{Walsh2} considered the wave equation in two space-time dimension:
\begin{equation}
\label{review9}
   \frac{\partial^{2} X }{\partial t^{2}}(t,x) -\frac{\partial^{2} X }{\partial x^{2}}(t,x) = \xi(t,x),
\end{equation}
whose solution is the Brownian sheet $w(s,t)$. For a Brownian motion $w_{t}$ we have the following (LIL):
\begin{equation}
\label{review10}
   \limsup_{h \rightarrow 0} \frac{w_{t+h}-w_{t}}{\sqrt{2h \log\log h}} = 1, \; \; {\rm a.e.}
\end{equation}
for each $t$ fixed. The extension of this estimate to the Brownian sheet is: for each $s$ fixed:
\begin{equation}
\label{review102}
   \forall t \geq 0 \; \; \limsup_{h \rightarrow 0} \frac{w_{s+h,t}-w_{s,t}}{\sqrt{2h \log\log h}} = \sqrt{t}, \; \;  {\rm a.e.}
\end{equation}
Then a singularity is defined when we take $s$ as a random variable $S(\omega)$ and when the modulus of continuity (or the LIL) given by (\ref{review10}) fails, 
i.e. when we have:
\begin{equation}
\label{review11}
    \limsup_{h \rightarrow 0} \frac{w_{S+h,t}-w_{S,t}}{\sqrt{2h \log\log h}} = \infty, \; \; a.e., 
\end{equation}
and the propagation of singularities result of Walsh \cite{Walsh1} states that that if (\ref{review11}) holds for some $(S,t_{0})$ 
(supposing that $S$ is $\sigma(w_{s,t}, t\leq t_{0})$-measurable) then it holds for all $(S,t)$ with $t\geq t_{0}$. This fact means that the singularities propagates in the $t$ direction
(vertically) and by symmetry also horizontally in the $s$ direction if we take points $(s, T)$ with a random time $T$. This result was extended by Carmona and Nualart \cite{CN} who considered a non-linear generalization of (\ref{review9}):
\begin{equation}
\label{review12}
   \frac{\partial^{2} X }{\partial t^{2}}(t,x) -\frac{\partial^{2} X }{\partial x^{2}}(t,x) = a(X(t,x)) \xi(t,x)+B(X(t,x)),
\end{equation}
and these authors studied the existence of the above mentioned kind of singularities for (\ref{review12}) and their propagation and reflexion with some other probabilistic properties. 
The results of Walsh were also extended by Blath and Martin \cite{BM} to semi-fractional Brownian sheets $X_{s,t}$ which are two-dimensional Gaussian random fields with a covariance of the form $E X_{s,t} X_{s',t'} = t\wedge t' (s^{\alpha}+s'^{\alpha} -|s-s'|^{\alpha})/2, \alpha \in ]0, 2[$.

\section{Random symbols and random pseudodifferential operators}

For the purpose of the next section, we introduce a class
of random symbols and random pseudodifferential operators. By a random 
symbol, we mean a symbol $p(x,\xi,\omega)$ which depends on a parameter $\omega\in
 \Omega$ and such that for each $x,\xi\in \R^{d}$ $\omega\mapsto p(x,\xi,\omega)$ is
measurable. 

We shall use the following classes of random symbols:
\\
\\
\noindent
$\bullet$ The class $RS^{m}$ is formed by the set of random symbols $p(x,\xi,\omega)$
such that for each compact $K\subset \R^{d}$ and $\alpha,\beta \geq 1$ there is a constant 
$C(\alpha,\beta, K, \omega)$ (measurable w.r.t $\omega$) that satisfies: 
\begin{equation}
\label{randomsymb1}
\forall x \in K, \xi \in \R^{d}: |D_{\xi}^{\alpha}D_{x}^{\beta}a(x,\xi, w)|\leq C(\alpha, \beta, K, \omega) (1+|\xi|)^{m-|\alpha|} \; \; a.e. 
\end{equation}
\noindent
We note here that the bounds of the derivatives of the symbols are allowed to be pathwise ($\omega$) dependent.
\\
\noindent
N.B. In the following, the letter $\omega$ may be omitted in the notation of random symbols.
\\
\\
\noindent
$\bullet$ The class $RS^{m}_{h}$ is formed by the set of random symbols $p(x,\xi,\omega) \in RS^{m}$ which are
homogeneous of degree $m$ in $\xi$ (at least for $\xi \geq 1$).
\\
\\
\noindent
$\bullet$ The class $RS^{m}_{Cl}$ that we call the classical random symbols (as in the deterministic case), 
is formed by the set of random symbols $p(x,\xi,\omega) \in RS^{m}$ such that there exist a sequence of random symbols 
$a_{j}(x,\xi,\omega) \in RS^{m-j}_{h}, j=0, 1, 2, ...$ with:
\begin{equation}
\label{RS1}
 p(x,\xi,\omega) - \sum_{j=0}^{n} a_{j}(x,\xi,\omega)  \in RS^{m - (n+1)} \; \; \forall N \geq 0, \; a.e. 
\end{equation}
When a random symbol $a(t, x,\xi)$ depends on the time
parameter $t$, the family $(a(t, x,\xi)), t\in \R$ will be assumed to be ${\cal F}_{t}$ adapted.
\\
As mentioned in the introduction, random symbols and random PDOs have already been used by Dedik, Shubin \cite{DS2}, Pankov \cite{Pankov}, Fedosov-Shubin \cite{FS}. 
These authors considered symbols $a(x,\xi, \omega)$ which depend on a random parameter $\omega$, satisfying the inequalities (\ref{randomsymb1}). 
However, the constants $C(\alpha, \beta, K, \omega)$ do not depend on $\omega$ in these works. 
On the other hand, Liu and Zhang \cite{LZ} introduced a class
of random PDOs and random symbols satisfying the inequalities (\ref{randomsymb1}) and depending on a time parameter $t\in [0,T]$ with constants $C(t, \alpha, \beta, K, \omega)$ that depend on $\omega$ and verify conditions like:
\begin{equation}
\label{RS2}
 \E \int_{0}^{T} |C(t, \alpha, \beta, K)|^{p} < \infty
\end{equation}

These authors transposed some of the usual properties of symbols and PDOs to this stochastic settings. In our case, we will be led to use symbols which are defined as a stochastic integral of random symbols: 
\begin{equation}
\label{RS3}
 q(t,x,\xi):= \int_{0}^{t}p(\t,x,\xi)dw(\t) 
\end{equation}

However, the condition (\ref{randomsymb1}) is not trivial for $q(t,x,\xi)$ as it involves a pathwise bound for stochastic integrals.
Nonetheless, this condition can be verified if we consider homogeneous random symbols; this is the subject of the following proposition:

\begin{prop}
\label{PropRS1}
Let $(p(t,x,\xi) \in RS^{m}_{h} , t\in [0,T])$ be a family of random homogeneous symbols such that  for each
$(x,\xi)$, $p(t,x,\xi)$ is a semimartingale and there exists a constant $\gamma > d$ with:
\\
(1) $\forall \alpha, \beta \geq 0, \; \exists C_{1} : \E |D_{\xi}^{\alpha}D_{x}^{\beta}p(t, x,\xi)|^{\gamma} \leq C_{1} < + \infty $
\\
(2) $\forall \alpha, \beta \geq 0 \; \exists C_{2}: \E |\< D_{\xi}^{\alpha}D_{x}^{\beta}p(t, x,\xi), w  \right>_{t} |^{\gamma} \leq C_{2} < + \infty $
\\
Then:
\[ q(t,x,\xi):= \int_{0}^{t}p(\t,x,\xi)\circ dw(\t) \]
(which is well defined) is a random symbol in $RS^{m}_{h}$.
\\
$\bullet$ The condition (2) is verified in particular if the martingale part of $p(t,x,\xi)$ is
of the form:
$$\int_{0}^{t}r(\t,x,\xi)dw(\t), \; {\rm  with} \;  r\in RS^{m}_{h} \; {\rm and}  \;
\E |D_{\xi}^{\alpha}D_{x}^{\beta}r(t, x,\xi)|^{\gamma} < + \infty $$
\end{prop}

NB. The brackets $\< U, V \right>_{t}$ in the point (2) of the proposition denotes here the quadratic variation of the processes $U,V$, see \cite{Ab1}. The proof of this proposition relies on the following lemmas:

\begin{lm}
\label{RSlemma1}
Let $g(x)$ be a random field indexed by $x\in \R^{d}$ and suppose that there exist constants $\gamma > d$ and $C>0$ such that:
\begin{equation}
\E | \frac{\partial g}{\partial x_{i}}(x)|^{\gamma} \leq C < +\infty \; \; \forall x\in K, i=1, \cdots, d
\end{equation}
where $K \subset \R^{d} $ is a compact set. Then $g$ has a modification (still denoted by $g$) which is a.e. continuous and for almost all $\omega\in \Omega$ there exists a constant $C(K, \omega)$ such that:
$$ \sup_{x\in K} |g(x)| \leq C(K, \omega) < + \infty. $$
\end{lm}
{\it Proof.} To prove the continuity of $g$ we use the multidimensional Kolmogorov-Centsov criterion (see Kallenberg \cite{Kallenberg}) that we recall here: if $(E,d)$ be a complete metric space, and $g(x), x\in K \subset \R^{d}$ an $E$-valued process such that there exist $C, \alpha, \beta > 0$:
\begin{equation}
\label{ContinuityCriterion}
\E (d(g(x), g(y))^{\alpha} \leq C |x - y|^{d+\beta} \; \; \forall x, y \in I,
\end{equation}
then $g(x)$ has a modification which is almost-surely $\lambda$-H\"older continuous for all $\lambda \in ]0, \beta/\alpha[$ and this modification verifies (\ref{ContinuityCriterion}).
\\
For the function $g$ of the lemma, we have:
\[ g(y)-g(x) = \int_{0}^{1}\sum_{k=1}^{d} \frac{\partial g}{\partial x_{k}}(x+t(y-x)) (y_{k}-x_{k})dt   \]
By the the H\"older inequality we have for $x_{k} \geq 0, p>1 $:
\begin{equation}
\label{Holder}
(\sum_{k=1}^{N} |x_{k}|)^{p}\leq N^{p-1} \sum_{k=1}^{N} |x_{k}|^{p}
\end{equation}
Then, for $\alpha \geq 1$:
\begin{align*}
|g(y)-g(x)|^{\alpha} \leq & 
      d^{\alpha-1} \sum_{k=1}^{d} (\int_{0}^{1} \frac{\partial g}{\partial x_{k}}(x+t(y-x))dt)^{\alpha} |y_{k}-x_{k}|^{\alpha} \\
			\leq & d^{\alpha-1} \sum_{k=1}^{d} \int_{0}^{1} |y_{k}-x_{k}|^{\alpha} (\frac{\partial g}{\partial x_{k}}(x+t(y-x)))^{\alpha} dt,
\end{align*}
which implies by the assumption of the lemma that:
\begin{equation*}
\E |g(y)-g(x)|^{\alpha} \leq d^{\alpha} C |y_{k}-x_{k}|^{\alpha}
\end{equation*}
If we choose $\alpha= \gamma> d, \beta = \gamma - d > 0$ then the criterion (\ref{ContinuityCriterion}) will be verified and there is a modification of $g$ which is a.e. continuous; this modification is therefore almost everywhere bounded on every compact $K$, which means that 
$$\exists C(K, \omega):  \sup_{x\in K} |g(x)| \leq C(K, \omega) < + \infty $$ 
for almost all $\omega$.$\Box$

\begin{lm}
\label{RSlemma2}
Let $f(\xi)$ be a random field indexed by $\xi\in \R^{d}$ and suppose that there exist constants $\gamma > d$ and $C>0$ such that:
\begin{equation}
\E | \frac{\partial f}{\partial x_{k}}(x)|^{\gamma} \leq C < +\infty \; \; \forall \xi\in U_{1}, k=1, \cdots, d
\end{equation}
where $U_{1}=\{\xi: |\xi|=1 \}$. Then $f$ has a modification (still denoted by $f$) which is a.e. continuous and for almost all $\omega\in \Omega$ there exists a constant $C(K, \omega)$ such that:
$$ \sup_{\xi \in \R^{d}} |f(\xi)| \leq C(\omega) |\xi|^{m} $$
\end{lm}
{\it Proof.} We make the change of variable $x=\xi/|\xi|$ and set $g(x)= f(x)=f(\xi/|\xi|)=f(\xi)/|\xi|^{m}$. 
Then $ \E |\partial_{k} g (x)|^{\gamma} = \E |\partial_{k} f (x)|^{\gamma}  \leq C $ for all $\in U_{1}$ and by Lemma \ref{RSlemma1}, there exists a modification of $g$ (and $f$), still denoted by the same letters, and a constant $C(\omega)$ such that $|g (x)| \leq C(\omega)$ for almost all $\omega$, i.e., $|f(\xi)| \leq C(\omega) |\xi|^{m}$ for all $\xi \in \R^{d}$. $\Box$
\\
\noindent
{\it Proof of Proposition \ref{PropRS1}.} First, we shall show that $q$ has a modification for which the 
derivatives $D^{\alpha}_{\xi}D^{\beta}_{x}q(t,x,\xi)$ exist and
\begin{equation}
\label{RS31}
D_{\xi}^{\alpha}D_{x}^{\beta}\int_{0}^{t}p(\t,x,\xi)\circ dw(\t)=
\int_{0}^{t}D_{\xi}^{\alpha}D_{x}^{\beta}p(\t,x,\xi)\circ dw(\t)
\end{equation}
Since this assertion is of local character (in $x$ and $\xi$) we may suppose that
all the derivatives $D_{\xi}^{\alpha}D_{x}^{\beta}p$ are bounded and
H\"older continuous. We can thus apply the theorem 1.2 of
Kunita \cite{K1} which gives (\ref{RS31}). 
We also have:
\[ D_{\xi}^{\alpha}D_{x}^{\beta}q(t,x,\xi)=\int_{0}^{t}
D_{\xi}^{\alpha}D_{x}^{\beta}p(\t,x,\xi) dw(\t) +\frac{1}{2}\int_{0}^{t}
D_{\xi}^{\alpha}D_{x}^{\beta}\left< p (\cdot,x,\xi),w (\cdot)  \right>_{\t}d\t, \]
which yields:
\begin{align}
\label{RS4}
 \E|D_{\xi}^{\alpha}D_{x}^{\beta}q(t,x,\xi)|^{\gamma} \leq &
2^{\gamma-1}[ \E |\int_{0}^{t} D_{\xi}^{\alpha}D_{x}^{\beta}p(\t,x,\xi) dw(\t)|^{\gamma} \\
       &  + \E|\frac{1}{2}\int_{0}^{t}D_{\xi}^{\alpha}D_{x}^{\beta}\left< p (\cdot,x,\xi),w (\cdot)  \right>_{\t}d\t|^{\gamma}],
\end{align}
where we have used the H\"older inequality ($\sum_{k=1}^{N} |x_{k}||y_{k}|\leq (\sum_{k=1}^{N} |x_{k}|^{p})^{1/p} (\sum_{k=1}^{N} |y_{k}|^{q})^{1/q}$ with $p=\gamma, N=2, y_{k}=1$). By the martingale moment inequalities and H\"older's inequality we get:
\begin{align*}
\E |\int_{0}^{t} D_{\xi}^{\alpha}D_{x}^{\beta}p(\t,x,\xi) dw(\t)|^{\gamma} 
         & \leq C_{\gamma} \E (\int_{0}^{t} |D_{\xi}^{\alpha}D_{x}^{\beta}p(\t,x,\xi)|^{2} d\t)^{\gamma/2} \\
         & \leq C_{\gamma} \E (\int_{0}^{T} |D_{\xi}^{\alpha}D_{x}^{\beta}p(\t,x,\xi)|^{\gamma} d\t) T^{\gamma/2-1} \\
\end{align*}
As for the second term of the r.h.s. of (\ref{RS4}), we use the H\"older inequality to get:
\[ \E|\frac{1}{2}\int_{0}^{t}
D_{\xi}^{\alpha}D_{x}^{\beta}\left< p (\cdot,x,\xi),w (\cdot)  \right>_{\t}d\t|^{\gamma}
\leq
\frac{1}{2^{\gamma}}  T^{\gamma-1}\int_{0}^{T} \E |D_{\xi}^{\alpha}D_{x}^{\beta}\left< p (\cdot,x,\xi),w (\cdot)  \right>_{\t}|^{\gamma}d\t
\]
Hence, by the assumptions (1) and (2) of the proposition, the last inequalities imply that for all $\alpha, \beta \geq 0 $ we have:
\[ \E|D_{\xi}^{\alpha}D_{x}^{\beta}q(t,x,\xi)|^{\gamma} \leq  T^{\gamma/2-1} T C_{\gamma} C_{1}  + T^{\gamma-1} T C_{2} \]
Now we use lemma \ref{RSlemma1} to control the derivatives w.r.t. $x \in K$ and Lemma \ref{RSlemma2} in order to control the derivatives w.r.t. $\xi\in \R^{d}$. For the later case, we also use the fact that if a differentiable function $f(\xi)$ is homogeneous of degree $m$, then the derivative $D^{\alpha}f(\xi)$ is homogeneous of degree $m-\alpha$. 
We will thus have: for any compact $K \subset \R^{d}$ and for all $\alpha, \beta \geq 0 $ and for almost all $\omega\in \Omega$, there exists a constant $C(\alpha, \beta, K, \omega)$ such that:
$$ |D_{\xi}^{\alpha}D_{x}^{\beta}p(\t,x,\xi)| \leq  C(\alpha, \beta, K, \omega) (1+|\xi|)^{m-|\alpha|}, $$
which means $q(t,x,\xi):= \int_{0}^{t}p(\t,x,\xi)\circ dw(\t)$ is a random symbol in $RS^{m}_{h}$.
As for the second part of the proposition, we note that $\left< p (\cdot,x,\xi),w (\cdot)  \right>_{\t}=r(\t,x,\xi)$ and condition (2) of the proposition is verified by the assumption on $r$. $\Box$
\\
\\
\noindent
A classical symbol can be constructed from a sequence of homogeneous symbols; we shall need the following extension of this result to
the random symbols: 
\begin{prop}
\label{propRS2}
Let $q_{j}(x,\xi,\omega) \in RS^{m-j}_{h}, j=0, 1, \cdots $ be a sequence of random symbols homogenous of degree $m-j$. 
Then there exists a random symbol $q(x,\xi,\omega) \in RS^{m}_{Cl}$ such that $q(x,\xi,\omega) \sim \sum_{j=0}^{\infty} q_{m-j}(x,\xi,\omega)$, a.e.
\end{prop}
Unlike the situation of Proposition \ref{PropRS1}, the proof in the deterministic case can be adapted to the case of random symbols without difficulties (see, e.g. \cite{SaintRaymond}). For completeness, this proof will be outlined in $\S 6$. 

\begin{lm}
\label{lemmaRS3}
Let $p(t,x,\xi)$ be a random symbol satisfying the assumptions of
Proposition \ref{PropRS1} and $v\in M^{2}(I,H^{s})$. Define $u(t)=\int_{0}^{t}v(\t)\circ
 dw(\t)$. Then
\begin{equation}
\label{fubini}
 q(t,x,D)u(t)=\int_{0}^{t}p(\t,x,D)u(\t)\circ dw(\t)+\int_{0}^{t}q(\t,x,D)
v(\t)\circ dw(\t).
\end{equation}
\end{lm}
{\it Proof.} Suppose first that $v(t)\in {\cal S}$ a.s. for all $t$. Then a
stochastic Fubini theorem (see, e.g., \cite{DZ}) yields:
\[ \hat{u}(t,\xi)=\int_{0}^{t}\hat{v}(\t,\xi)\circ dw(\t). \]
where $\hat{u}(t, \cdot)$ means the Fourier transform of $u(t, \cdot)$ (w.r.t. $x$), and
\[
q(t,x,D)u(t)=\int_{\R^{d}}q(\t,x,\xi)[\int_{0}^{\t}\hat{v}(\theta,\xi)\circ 
            dw(\theta)]e^{i x.\xi}d\xi. \]
But the It\^o formula gives
\[ q(\t,x,\xi)\int_{0}^{\t}\hat{v}(\theta,\xi)\circ dw(\theta)=\int_{0}^{\t}
  q(\t,x,\xi)\hat{v}(\theta,\xi)\circ dw(\theta)+\int_{0}^{\t}
  p(\t,x,\xi)\hat{u}(\theta,\xi)\circ dw(\theta).\]
Using again the stochastic Fubini theorem we deduce that (\ref{fubini}) holds
for $v(t)\in {\cal S}$. The case where $v\in M^{2}(I, H^{s})$ follows by
a density argument. $\Box$

The next proposition concerns the dependance of a particular stochastic flow with respect to the initial conditions:
\begin{prop}
\label{propRS3}
Let $a(t,x,\xi), t\in [0,T]$ be a bounded family in $S^{1}$ and
consider the stochastic  equations 
\[
({\cal C}):
\left\{ \begin{array}{c}
 \displaystyle dx(t)=\frac{\partial a}{\partial \xi}(t,x(t),\xi(t))\circ dw(t) \\
\displaystyle d\xi(t)=-\frac{\partial a}{\partial x}(t,x(t),\xi(t))\circ dw(t)
\end{array}\right.
\]
with $x(0)=x_{0},\xi(0)=\xi_{0}$. Then $({\cal C})$ has a global solution 
defined on $[0,T]$. Furthermore, if we denote  $\phi_{t}(x_{0},\xi_{0})=
(x(t),\xi(t))$ the stochastic flow of diffeomorphisms associated to 
(${\cal C}$) and
$(\bar{x}(t),\bar{\xi}(t))=\phi^{-1}_{t}(x_{0},\xi_{0})$, then we have for $n\geq 1$:
\\
$\bullet$ The quantities $\E |x(t,x_{0},\xi_{0})|^{n}, \E |\xi(t,x_{0},\xi_{0})|^{n}, \E |\bar{x}(t,x_{0},\xi_{0})|^{n}, \E |\bar{\xi}(t,x_{0},\xi_{0})|^{n}$ are bounded.
\\
$\bullet$ For all $\alpha, \beta \geq 1$, the following quantities are bounded:
\[ E |D_{\xi_{0}}^{\alpha}D_{x_{0}}^{\beta} x(t,x_{0},\xi_{0})|^{n}, 
E |D_{\xi_{0}}^{\alpha}D_{x_{0}}^{\beta}\bar{x}(t,x_{0},\xi_{0})|^{n}, 
E |D_{\xi_{0}}^{\alpha}D_{x_{0}}^{\beta} \xi (t,x_{0},\xi_{0})|^{n},  
E |D_{\xi_{0}}^{\alpha}D_{x_{0}}^{\beta}\bar{\xi}(t,x_{0},\xi_{0})|^{n},
\]
$\bullet$ If $a(t,x,\xi)$ is homogeneous of degree $1$ in $\xi$, then $\phi_{t}(x_{0},\xi_{0})$ is homogeneous of degree $1$ in $\xi_{0}$
\end{prop}
This proposition and the following corollary will be used in the proof of the main result. Their proofs are given in the section 6.
\begin{co}
\label{coroRS1}
Let $p(x,\xi)\in S^{m}_{h}$ and $\phi_{t}$ be the flow associated to (${\cal C}$).
 Then
$q_{t}(x,\xi):= p(\phi_{t}(x,\xi)), l_{t}(x,\xi)=p(\phi_{t}^{-1}(x,\xi))$ define a family of random symbols in $RS^{m}_{h}$.
\end{co}
\section{Propagation of singularities for hyperbolic stochastic partial differential equations}

In this section we study the singularities of the solutions to Eq. (E), at the wave front level. We prove a propagation
of singularities result which is similar to the deterministic case. For the sake of simplification we state and prove this result in the scalar case ($d'=1$).
\\
\noindent
Let us consider the stochastic equation:
\begin{equation}
\label{hyperbolicSPDE}
du(t)= A(t,x,D)u(t)\circ dw(t)+ B (t,x,D)u(t)dt , \; u(0)=u_{0} \in (H^{s}(\R^{d})).
\end{equation}
We set $A(t,x,\xi) = i a(t,x,\xi), B(t,x,\xi) = i b(t,x,\xi)$ and we suppose in this section that $a(t,x,\xi), b(t,x,\xi)$ satisfy the following additional condition:

$(v)$ The principal symbols of $a(t,x,\xi)$ and $b(t,x,\xi)$ are real and 
\[ a(t,x,\xi)\sim \sum_{j=0}^{\infty}a_{1-j}(t,x,\xi),\; 
b(t,x,\xi)\sim \sum_{j=0}^{\infty}b_{1-j}(t,x,\xi),\]
with $a_{1-j}(t,x,\xi), b_{1-j}(t,x,\xi)$ homogeneous in $\xi$ of degree $1-j$.
We are interested in the relationship between the singularities of $u(t, \cdot)$ 
and 
those of $u_{0}$. In the deterministic case ($a=0$) the singularities of the solution of the hyperbolic PDE: $\partial_{t}u + b(t,x,D)u=0$
$u(t, \cdot)$ propagate along the bi-characteristic curves of the principal symbol $b_{1}$ of $b$, i.e. the integral
curves of the Hamilton vector field 
\[ H_{b_{1}}=\frac{\partial b_{1}}{\partial \xi_{i}}\frac{
\partial}{\partial x^{i}}- \frac{\partial b_{1}}{\partial x^{i}}\frac{
\partial}{\partial \xi_{i}}. \]
In other words, if $\chi_{t}(x,\xi)$ represents the integral curves (or flow) associated to  $H_{b_{1}}$, then we have $WF(u(t, \cdot))=\chi_{t} (WF(u(0, \cdot))$. 
This theorem has been proved in several ways, see, e.g., Duistermaat-H\"ormander \cite{DH}, H\"ormander\cite{H2}-\cite{H5}, Sogge \cite{Sogge}, Rauch\cite{Ra}, Raltson\cite{Raltson}, Taylor\cite{Taylor1}. 

In the case of hyperbolic SPDEs (\ref{hyperbolicSPDE}), we have a similar result. We shall use random symbols and random PDOs to adapt the proof given in H\"ormander\cite{H5}, Raltson\cite{Raltson}, and it is instructive to recall the idea of this proof: given a couple $(x_{0}, \xi_{0}) \notin WF(u_{0})$, we want to prove that $(x_{t}, \xi_{t}):=(x, \xi):=(\chi_{t}(x_{0}, \xi_{0})) \notin WF(u_{t,\cdot})$, where $(x_{t}, \xi_{t})$ is the integral curve of the vector field $H_{b_{1}}$ (or Hamiltonian flow). We seek a family of operators $Q(t,x,D)$ which nearly commute with $iD_{t}u + b(t,x,D)u$, so that:
\begin{equation}
\label{propag1}
(iD_{t} + b(t,x,D)) Q(t,x,D) u - Q(t,x,D) (iD_{t} + b(t,x,D)) u \in S^{-\infty} 
\end{equation}
If $u(t,x)$ is a solution of $(iD_{t} + b(t,x,D)) u=0$, then we will have 
\begin{equation}
\label{propag2}
(iD_{t} + b(t,x,D)) Q(t,x,D) u \in S^{-\infty}
\end{equation}
This implies that $Q(t,x,D) u (t,x) \in C^{\infty}$. Now if we find $Q(t,x,D)$ as operators with classical symbols 
$Q (t, x,\xi) \sim \sum_{0}^{\infty}Q_{j}(t, x,\xi)$ then we can conclude that $(x, \xi) \notin WF(u_{t,\cdot})$ provided that
$Q_{0}(t, x,\xi) \neq 0$. The point is that when such $Q_{j}$ are constructed, it turns out that they verify: 
$$ Q_{0}(t, x_{t}, \xi_{t}):=Q_{0}(t, \chi_{t}(x_{0}, \xi_{0}))= q_{0} (x_{0}, \xi_{0}).$$
Since $q_{0}(x_{0}, \xi_{0}) \neq 0$, this means that we will have $ Q_{0}(t, x, \xi) \neq 0$ as soon as we have $ (x, \xi)=\chi_{t}(x_{0}, \xi_{0})$, i.e, when $(x, \xi)$ belongs to the integral curves of $H_{b_{1}}$. This proves that $WF(u_{t,\cdot}) \subset \chi_{t}(WF(u_{0,\cdot}))$; the equality of the two sets is obtained by time reversal. 

\begin{thm}
Let $u_{0}\in H^{s}$ and $u\in M^{s}_{2}(I,H^{s})$ be the solution to 
(\ref{hyperbolicSPDE}). Then $WF(u(t, \cdot))=\Phi_{t}(WF(u_{0}))$ a.e. where the transformation 
$\Phi_{t}(x_{0},\xi_{0})=(x_{t},\xi_{t})$ is given by
\[ dx(t)=\frac{\partial a_{1}}{\partial\xi_{i}}(t,x(t),\xi(t))\frac{\partial}{
\partial x^{i}}
\circ dw(t) +\frac{\partial b_{1}}{\partial\xi_{i}}(t,x(t),\xi(t))
\frac{\partial}{\partial x^{i}}
dt ,\]
\[d\xi(t)=-\frac{\partial a_{1}}{\partial x^{i}}(t,x(t),\xi(t))\frac{\partial}{
\partial \xi_{i}}
\circ dw(t) -\frac{\partial b_{1}}{\partial x^{i}}(t,x(t),\xi(t))
\frac{\partial}{\partial \xi_{i}}dt
\]
with $x(0)=x_{0}, \xi(0)=\xi_{0}$ and $a_{1}, b_{1}$ are the principal symbols of order $1$ of $a$ and $b$.
\end{thm}
{\it Proof.}
To simplify the notations we suppose that $b\equiv 0$, as we are mainly interested in the proof related to the stochastic part. 
Let $(x_{0}, \xi_{0})\notin WF(u_{0})$. Then there is a symbol
$q (x,\xi) \sim \sum_{0}^{\infty}q_{j}(x,\xi)\in S^{0}_{ph}$ with $q_{0}(x_{0},\xi_{0})\neq
 0$ and $q(x,D)u_{0}\in C_{0}^{\infty}$. 
\\
To prove the theorem we shall construct a family of random symbols such that $Q (t,x,\xi) \sim \sum_{0}^{\infty}Q_{j}(t,x,\xi)\in S^{0}_{ph}$ a.e.
with $Q(t, x,D)u_{t}\in C_{0}^{\infty}$ and $Q_{0}(t, x_{t},\xi_{t}) = Q_{0}(t, \Phi_{t}(x_{0},\xi_{0}))\neq 0$. 

Let $p(t,x,\xi) \sim \sum_{0}^{\infty}p_{j}(t,x,\xi)\in S^{0}_{ph}$ be a random symbol in $RS^{0}$ which satisfies the 
conditions of Proposition \ref{PropRS1}; we shall seek $Q$ of the form: 
\[ Q(t,x,\xi)= q(x,\xi)+\int_{0}^{t}p(\t,x,\xi)\circ dw(\t). \] 
Then by Lemma \ref{lemmaRS3}  we have:
\begin{align*}
Q(t) u(t)-Q(0)u_{0} =&\int_{0}^{t}p(\t)u(\t)\circ dw(\t)
          +\int_{0}^{t}Q(\t)a(\t)u(\t)\circ dw(\t) \\
   =& \int_{0}^{t}a(\t)Q(\t)u(\t)\circ dw(\t)\\
    & +\int_{0}^{t}(Q(\t)p(\t) +[Q,a](\t))u(\t)\circ dw(\t). 
\end{align*}
We shall choose $p$ such that $c(t,x,\xi):= p(t,x,\xi)+[Q,a](t,x,\xi)\in S^{-\infty}$ a.e. We have:
\begin{align}
\label{commut}
[Q,a](t,x,\xi) & \sim \sum_{|\alpha | \geq 1} \frac{i^{\alpha}}{|\alpha| !} (D^{\alpha}_{\xi} Q D^{\alpha}_{x} a  
                            - D^{\alpha}_{\xi} a D^{\alpha}_{x} Q ) \\
               & \sim \{a_{1}, Q\} + \{a- a_{1}, Q\} + \sum_{|\alpha | \geq 2} \frac{i^{\alpha}}{|\alpha| !} (D^{\alpha}_{\xi} a D^{\alpha}_{x} Q  - D^{\alpha}_{\xi} Q D^{\alpha}_{x} a ) \nonumber
\end{align}
In this sum, the term of order $0$ is $\{Q, a_{1}\} = H_{a_{1}} Q_{0}(t,x, \xi)$, where $H_{a_{1}}= (\partial a_{1}/\partial x)\partial/\partial \xi- (\partial a_{1}/\partial \xi)\partial/\partial x$ and the principal symbol of $c(t,x,\xi)$ (order $0$) is:
\[ c_{0}(t,x,\xi)= p_{0}(t,x,\xi)+ H_{a_{1}} Q_{0}(t,x,\xi),\]
 
The symbol of order $-j$ of $[Q,a](t,x,\xi)$ contains the terms $D^{\alpha}_{\xi} Q D^{\alpha}_{x} a - D^{\alpha}_{\xi} a D^{\alpha}_{x} Q$ with $|\alpha|\leq j+1$, and when taking the terms $Q_{j'}, a_{j''}$ of the expansion of $Q$ and $a$, the terms that contribute to the order $-j$ of $[Q,a](t,x,\xi)$ should satisfy: $j', j'' \leq j-1$. 
Hence the symbol of order $-j, j\geq 1$ of $c(t,x,\xi)$ can be written as:
\[  c_{j}(t,x,\xi)= p_{j}(t,x,\xi)+ H_{a_{0}} Q_{j}(t,x,\xi) +R_{j}(t,x,\xi),\]
where the symbols $R_{j}(t,x,\xi)$ are determined by $Q_{0},..., Q_{j-1}$, this can be seen by expanding (\ref{commut}), see Raltson \cite{Raltson} p. 96 for explicit expressions of similar quantities.
\\
Let us remark that, in order to have $c(t,x,\xi) \in S^{-\infty}$ a.e., the condition $c_{j}(t,x,\xi)=0$ for all $t$ is equivalent to
$\int_{0}^{t}c_{j}(\t,x,\xi)\circ dw(\t)=0$ a.e. for all $t$. Thus, we have
to determine the $(Q_{j})$ which satisfy
\begin{equation}
\label{hj1}
dQ_{0}(t,x,\xi)= -H_{a_{1}}Q_{0}(t,x,\xi)\circ dw(t) ,
\end{equation}
\begin{equation}
\label{hj2}
dQ_{j}(t,x,\xi)= -H_{a_{1}}Q_{j}(t,x,\xi)\circ dw(t) +R_{j}(t,x,\xi)\circ dw(t).
\end{equation}

The stochastic characteristic equations associated to (\ref{hj1}) are the
Eqs. $({\cal C})$, which admit a global solution by Proposition \ref{propRS3}. Hence,
by Kunita \cite{K1}, the equation (\ref{hj1}) has a global solution which is
 given by $Q_{0}(t,x,\xi)=q_{0}(\Phi_{t}^{-1}(x,\xi))$. By the same arguments,
(\ref{hj2}) has a global solution given by:

\begin{equation}
\label{hj22}
 Q_{j}(t,x,\xi)=Q_{j}((\Phi_{t}^{-1}(x,\xi))+\int_{0}^{t}R_{j}(\t, \Phi_{\t}\circ \Phi_{t}^{-1}(x,\xi))\circ dw(\t).
\end{equation}

Now let us check the properties of the random symbols $Q_{j}(t,x,\xi)$ constructed by the above formula. First, observe that
the random maps $\Phi_{t}(x,\xi)$ and $\Phi_{t}^{-1}(x,\xi)$ are homogeneous of degree 0 and 1 respectively in $x$ and $\xi$. 
On the other hand an inspection of the terms involved in the random symbols $R_{j}(t, x,\xi)$ shows that these symbols are homogeneous of degree $-j$ in $\xi$.
By corollary \ref{coroRS1}, the first terms of the r.h.s. of (\ref{hj22}), $Q_{j}((\Phi_{t}^{-1}(x,\xi))$.
form a set of random symbols in $RS^{-j}_{h}$. As for the second term $\int_{0}^{t}R_{j}(\t, \Phi_{\t}\circ \Phi_{t}^{-1}(x,\xi))\circ dw(\t)$, 
since $R_{-j}(t, x,\xi) \in RS^{-j}_{h} $ as we have just seen, we can apply Proposition \ref{PropRS1}; its conditions are indeed verified: 
for Condition (1) we have to verify that there is a constant $\gamma > d$ such that:
\begin{equation}
\label{condition1}
 \forall \alpha, \beta \geq 0, \exists C_{1} : \E |D_{\xi}^{\alpha}D_{x}^{\beta}R_{j}(t, \Phi_{\t}\circ \Phi_{t}^{-1}(x,\xi))|^{\gamma} \leq C_{1} < + \infty
\end{equation}
To see this, we remark that (\ref{condition1}) is verified with $\gamma= d+1$ when $R_{j}= Q_{0}$; by Proposition \ref{propRS2} and Corollary \ref{coroRS1}, we see that
this condition holds for  $R_{1}$ which depends on $Q_{0}$, and the same argument can be applied by induction to show that (\ref{condition1}) is verified for $R_{j}$ 
which depends on $Q_{j'}, j'=0, ..., j-1$.
\\
As for the condition (2) of Proposition \ref{PropRS1}, we use the It\^o formula to expand $R_{j}(t, \Phi_{\t}\circ \Phi_{t}^{-1}(x,\xi))$ and we note that 
its martingale part  is of the form $\int_{0}^{t}r_{j}(\t,x,\xi)dw(\t)$ with $r_{j}\in RS^{-j}_{h}$ and by the same arguments used for condition (\ref{condition1})
we will have $\E |D_{\xi}^{\alpha}D_{x}^{\beta}r_{j}(t, x,\xi)|^{d+1} < + \infty $.
\\
Now since $Q_{j}(t,x,\xi)\in RS^{-j}_{h}$, by Proposition \ref{propRS2} there exists a random symbol $Q(t,x,\xi) \in RS^{0}_{Cl}$ such that 
$Q(t,x,\xi)\sim \sum_{j=0}^{\infty} Q_{j}(t,x,\xi)$. 
With this choice of $Q$ we have:
\[ Q(t)u(t)-Q(0)u_{0}=\int_{0}^{t}a(\t)Q(\t)u(\t)\circ dw(\t)+\int_{0}^{t} c(\t)u(\t)\circ dw(\t),\]
with $c(\t)\in S^{-\infty}$ for all $\t$, which implies that $q(t,x,D)u(t)\in 
H^{s}$ a.e. for all $s$ by Theorem 2.1 of \cite{Ab1}. To summarize, given $(x_{0},\xi_{0})
\notin WF(u_{0})$ and $q(x,\xi)\in S^{0}_{ph}$ with $q_{0}(x_{0},\xi_{0})\neq
0$ and $q(x,D)u_{0}\in C_{0}^{\infty}$, we have constructed a symbol 
$q(t,x,\xi)\in RS^{0}_{ph}$ such that $q(t,x,D)u(t)\in C_{0}^{\infty}$ and
$q_{0}(t, \Phi_{t}(x_{0},\xi_{0}))=q_{0}(x_{0},\xi_{0})\neq 0 $ which
implies that $\Phi_{t}(x_{0},\xi_{0})\notin WF(u(t, \cdot))$ and $WF(u(t, \cdot))
\subset \Phi_{t}(WF(u_{0}))$. To prove the converse, let $t\in ]0,T]$ be
 fixed and let us denote by $U_{b}(t,s)\phi, s\in [0,t]$ the solution of the backward 
equation
\[ u(s)=\phi-\int_{s}^{t}a_{\t}(x,D)u(\t)\circ \hat{d}w(\t). \]
By proposition 4.5 of \cite{Ab1} we have $U_{b}(t,0)u(t)=u_{0}$ a.s. On the other hand,
given a random symbol $q(t)=q(t,x,\xi)$ one 
can construct a family of random symbols $\hat{q}(\t), 0\leq \t\leq t$ such 
that 
\[ \hat{q}(t)\phi-\hat{q}(s)u(s)-\int_{s}^{t}a(\t)\hat{q}(\t)u(\t)\circ 
\hat{d}w(\t)=\int_{s}^{t}c(\t)u(\t)\circ \hat{d}w(\t),\]
with $c(\t)\in S^{-\infty}$. This can be done as above by using backward
equations. From this we conclude in the same way that for each (deterministic)
$\phi$ in some $H^{s}$ we have $WF(U(t,0)\phi)\subset \Phi_{-t}(WF(\phi))$
 for almost all $\omega$. Now let $\omega$ be given and fix
$\phi=u(t,\omega)$ in some $H^{s}$, then we have $WF(u_{0})=WF(
U(t,0,\omega)\phi)\subset \Phi_{-t}WF((u(t,\omega)))$ i.e. $\Phi_{t}(WF(u_{0}))
\subset WF(u(t,\omega))$. $\Box$
\\
\\
\noindent
{\it Remark: The case of differential operators.} We shall consider the case
 where $a_{t}(x,D), b_{t}(x,D)$
are differential operators to obtain simply a 'majorization' of the wave
front set of the solution to $({\cal E})$ as in the deterministic case. For
notational convenience we only consider the equation
\[ du(t)= \alpha(t,x)\frac{\partial u}{\partial x}\circ dw(t)+\beta(t,x)
\frac{\partial u}{\partial x}dt,\; \; u(0)=u_{0},\]
whose solution is $u(t,x)=u_{0}(\phi_{t}^{-1}(x))$ where $\phi_{t}(x)$ is the
stochastic flow associated to the eq. $dX(t)=\alpha(t,X(t))\circ dw(t)+
\beta(t,X(t)dt$.

\noindent
Then for each $w$, $u(t,x)$ can be written as a Fourier integral distribution
\[ u(t,x)=K_{t}u_{0}(x):= \left< K_{t}(x, \cdot), u_{0} (\cdot)   \right>, \]
with
\[ K_{t}(x,y)= (2\pi)^{-n}\int e^{i(\phi_{t}^{-1}(x)-y)\xi}d\xi .\]
It is known that the singularities of $K_{t} u_{0}$ are given by
\begin{eqnarray*}
WF(K_{t}u_{0})&=& WF'(K_{t})\circ WF(u_{0})\\
              &:=& \{(x,\xi): \exists (y,\eta)\in WF(u_{0}): (x,y,\xi,\eta)
                \in WF(K_{t})\}.
\end{eqnarray*}
On the other hand the wave front set of an oscillatory integral of the 
form $I(x)=\int a(x,\xi)\exp(i\phi(x,\xi))d\xi$ satisfies $WF(I)\subset
\{ (y,\eta): \partial \phi/\partial \xi(y,\eta)=0,\; \eta=\partial 
\phi/\partial x(y,\eta) \} $. Hence: 
\[ WF(K_{t})\subset \{ (x,\phi_{t}^{-1}(x),\partial_{x} \phi_{t}^{-1}. \xi,
-\xi), \; (x,\xi)\in (\R^{d})^{2}\}. \]
From this, it follows easily that $WF(K_{t} u_{0})\subset \{(x(t),\xi(t))=
(\phi_{t}(x),\partial_{x} \phi_{t}^{-1}(x)\xi):\; (x,\xi)\in WF(u_{0}) \}.$
Now we can verify that the above curves $(x(t),\xi(t))$ are the
solutions to the stochastic bicharateristic equations of the theorem. Indeed
we have:
\[ dx(t)=d\phi_{t}(x)=dx(t)=-\frac{\partial a_{1}}{\partial \xi}
(t,x(t),\xi(t))\circ dw(t) - \frac{\partial b_{1}}{\partial \xi}
(t,x(t),\xi(t))dt \]
since in this case $ia(t,x,\xi)=-i\xi\alpha(t,x),\; ib(t,x,\xi)=-i\xi
\beta(t,x)$. On the other hand, $Z_{t}:=\partial_{x} \phi_{t}(x)^{-1}$ 
satisfies (see, e.g. \cite{K2}, \cite{IW}):
\[dZ(t)=Z(t)\partial_{x}\alpha(t,x(t))\circ dw(t)-
  Z(t)\partial_{x}\beta(t,x(t))dt, \]
which implies that $\xi(t)= \partial_{x} \phi_{t}(x)^{-1}\eta$ verifies
the bicharacteristic equation of the theorem.

\section{Proofs of technical results }

{\it Proof Proposition \ref{propRS2}.}\\
\noindent
Let $\chi \in C^{\infty}(\R)$ such that: $\chi(t)=1$ if $|t|\leq 1$ and $\chi(t)=0$ if $|t|\geq 2$. We introduce $ \tilde{q}_{j}(x,\xi,\omega)= (1-\chi(\epsilon_{j}|\xi|)) q_{j}(x,\xi,\omega)$, with $\epsilon_{j}$ a decreasing sequence whose limit is $0$, that will be determined later (it will be $\omega$-dependent), and we set:
\[  q(x,\xi,\omega) = \sum_{j=0}^{\infty} \tilde{q}_{m-j}(x,\xi,\omega) \]
Then $q(.,.,\omega) \in C^{\infty}(\R^{d}\times \R^{d}) $ because for each $\xi$ the previous sum is finite, and the fact that $\tilde{q}_{j}- q_{j}=0$ when $|\xi| > 1/2\epsilon_{j}$ implies that $\tilde{q}_{j}- q_{j} \in RS^{-\infty} $ and $\tilde{q}_{j} \in  RS^{m-j}$ like ${q}_{j}$. To control the derivatives of $\tilde{q}$ we begin with those of $\tilde{q}_{j}$ and we suppose in the following that $|\xi| \geq 1$; we have:
\[ |\partial_{\xi}^{\alpha}\partial_{x}^{\beta}\tilde{q}_{j}| \leq \sum_{\gamma\leq\alpha}\frac{\alpha!}{\gamma!(\alpha-\gamma)!}\epsilon_{j}^{\gamma} 
                                                                      |\partial_{\xi}^{\alpha-\gamma}\partial_{x}^{\beta} q_{j}|
																										\leq C(\alpha,\beta, j, \omega)|\xi|^{m-j-|\alpha|}.
\]
To obtain the last inequality we use the following facts: $|\partial_{\xi}^{\alpha-\gamma}\partial_{x}^{\beta} q_{j}| \ C_{j}(\alpha, \beta, \omega) |\xi|^{m-j-|\alpha-\gamma|} \leq C_{j}(\alpha, \beta, \omega) |\xi|^{m-j-|\alpha|+|\gamma|} $ and for $ |\xi|\leq 2/\epsilon_{j} $ we have $|\xi|^{\gamma} \epsilon_{j}^{\gamma}\leq 2^{\gamma}$ while for 
$ |\xi| > 2/ \epsilon_{j} $ we have $\tilde{q}_{j}= q_{j}$.
\\
On the other hand if $\epsilon_{j}|\xi| <1 $ then $\tilde{q}_{j}= 0$; hence we may suppose that $\epsilon_{j}|\xi| \geq 1 $ as far as we consider $\tilde{q}_{j}$ and its derivatives and we have:
\[
 |\partial_{\xi}^{\alpha}\partial_{x}^{\beta}\tilde{q}_{j}| \leq  C(\alpha,\beta, j, \omega) \epsilon_{j}|\xi| |\xi|^{m-j-|\alpha|}
\]
and with the choice of an $\epsilon_{j}(\omega)= \min \{1/ C(\alpha,\beta, j, \omega), \; {\rm with} \; |\alpha+\beta|\leq j \}$ we have: $ |\partial_{\xi}^{\alpha}\partial_{x}^{\beta}\tilde{q}_{j}| \leq |\xi|^{1-j} |\xi|^{|\alpha|-m}  $ for $|\alpha+\beta|\leq j$. 
Now for $k\geq 0 $, to prove that $q- \sum_{j=0}^{k-1}q_{j}\in RS^{m-k}$ we set:
$q- \sum_{j=0}^{k-1}q_{j} = r^{(1)}_{k}+r^{(2)}_{k}+r^{(3)}_{k}$ with $r^{(1)}_{k}=\sum_{j=0}^{k-1}(\tilde{q}_{j})-q_{j})$, 
$ r^{(2)}_{k}=\sum_{j=k}^{N-1}\tilde{q}_{j}$ and $ r^{(3)}_{k}=\sum_{j=N}^{\infty}\tilde{q}_{j}$, where $N=\max(|\alpha+\beta|, k+1)$. We have $r^{(1)}_{k}\in RS^{-\infty}$ and $r^{(2)}_{k}\in RS^{m-k}$ as a finite sum of terms in $RS^{m-k}$. For the last term $r^{(3)}_{k}$, since $j\geq N\geq= \max(|\alpha+\beta|, k+1)$ in its sum, we have: 
\[
|\partial_{\xi}^{\alpha}\partial_{x}^{\beta}r^{(3)}_{k}|\leq \sum_{j=k+1}^{\infty} |\xi|^{|\alpha|-m}  |\xi|^{1-j} 
                                      \leq  |\xi|^{|\alpha|-m-k} \sum_{j=0}^{\infty} |\xi|^{-j}
\]
Hence for $|\xi|\leq 2$ we have $|\partial_{\xi}^{\alpha}\partial_{x}^{\beta}r^{(3)}_{k} \leq |\xi|^{|\alpha|-m-k} \sum_{j=0}^{\infty} 2^{-j}$. This shows that 
$q(x,\xi,\omega) - \sum_{j=0}^{k-1}q_{j} (x,\xi,\omega) \in RS^{m-k}$ a.e. $\Box$
\\
\\
\noindent
{\it Proof of Proposition \ref{propRS3}.}
\\
\noindent
$\bullet$ Let $X(t)=(x(t),\xi(t))$. If we write $({\cal C})$ in the form
$dX(t)=f(X(t))\circ dw(t)$, then from the assumptions on $a^{1}$, we deduce that $f$ is locally Lipshitz and has at most a linear growth as 
$x,\xi\rightarrow \infty$. By standard results on stochastic differential equations (cf., e.g., \cite{IW}), (${\cal C}$) has a global solution.

\noindent
$\bullet$ First, we write $x(t):=x(t,x_{0},\xi_{0}), \xi(t)=\xi(t,x_{0},\xi_{0})$ in the It\^o form:
\begin{align*}
dx(t)=& D_{\xi} a(t) dw(t) + \frac{1}{2}[D_{x}D_{\xi} a(t)D_{\xi} a(t) - D_{\xi}^{2} a(t)D_{x} a(t) ] dt \\
d\xi(t)=&-D_{x} a(t) dw(t) - \frac{1}{2}[D_{x}^{2} a(t)D_{\xi} a(t) - D_{x}D_{\xi} a(t)D_{x} a(t)] dt
\end{align*}
By the It\^o formula we have:
\begin{align*}
|x(t)|^{n}=& |x_{0}|^{n}+  \int_{0}^{t} n |x(\t)|^{n-1} D_{\xi} a(\t) dw(\t) + \frac{1}{2}\int_{0}^{t} n(n-1) |x(\t)|^{n-2} |x(\t)|^{n-2} (D_{\xi} a(\t))^{2} d\t \\
          & + \frac{1}{2}\int_{0}^{t} n |x(\t)|^{n-1} [D_{x}D_{\xi} a(\t)D_{\xi} a(\t) - D_{\xi}^{2} a(\t)D_{x} a(\t) ] d\t.
\end{align*}
We have a similar formula for $|\xi(t)|^{n}$ Let us set $\tilde{x}(t):=E\sup_{\t\leq t}|x(\t)|^{2n}$ and $\tilde{\xi}(t):=E\sup_{\t\leq t}|\xi(\t)|^{2n}$. Using the fact that
$ |x(\t)|^{n-1}$ and $|x(\t)|^{n-2}$ are  $\leq 1+|x(\t)|^{n}$, and the fact that $|D_{\xi} a(\t)|, |D_{x}D_{\xi} a(\t)|$ are bounded and applying martingale and H\"older inequalities we get:
\[ \tilde{x}(t) \leq |x_{0}|^{n}+ C_{2} \int_{0}^{t} \tilde{x}(\t) d\t, \]
which yields $ \tilde{x}(t)  \leq |x_{0}|^{n} e^{C_{2} T}$ by the Gronwall lemma. Hence $E\sup_{\t\leq t}|x(\t)|^{2n}$ is bounded. We use this result and the same previous argument to get the boundedness of $E\sup_{\t\leq t}|\xi(\t)|^{2n}$.

$\bullet$ Now, we shall sketch the proof of the estimates concerning 
$E |D_{\xi_{0}}^{\alpha}D_{x_{0}}^{\beta} x(t,x_{0},\xi_{0})|$ and 
$E |D_{\xi_{0}}^{\alpha}D_{x_{0}}^{\beta} x(t,x_{0},\xi_{0})|$. This will be done
by induction on $\alpha,\beta$. For simplicity we consider only
the quantities $E |D_{\xi_{0}}^{\alpha} x(t,x_{0},\xi_{0})|,
E |D_{\xi_{0}}^{\alpha}x(t,x_{0},\xi_{0})|$ (the proof is similar if we take into 
account the derivatives w.r.t. $x$). We use the notation:
\[ a(t) = a(t,x(t),\xi(t))\]
\noindent
For $\alpha=1$, we use Theorem 1.1 in Kunita \cite{K1} to get:
\begin{eqnarray*}
D_{\xi_{0}}x(t)&=& D_{\xi_{0}}x(0) + \int_{0}^{t}[D_{x}D_{\xi}a(\t) D_{\xi_{0}}x(\t)+D_{\xi}^{2}a(\t)D_{\xi_{0}}\xi(\t)]dw(\t)\\
   &+&\frac{1}{2}\int_{0}^{t}[ D_{\xi}D_{x}^{2}a(\t) D_{\xi}a(\t)+(D_{\xi} D_{x}a(\t))^{2}\\
   &-& D^{2}_{\xi}D_{x}a(\t)D_{x}a(\t)-D^{2}_{\xi} a(\t)D_{x}^{2} a(\t)]D_{\xi_{0}} x(\t) d\t\\
   &+&\frac{1}{2} \int_{0}^{t}[D_{\xi}^{2}D_{x}a(\t) D_{\xi}a(\t)+D_{\xi}D_{x} a(\t) D_{\xi}^{2}a(\t)\\
  &-&D^{3}_{\xi}a(\t)D_{x}a(\t)-D^{2}_{\xi} a(\t)D_{x}D_{\xi}a(\t)]D_{\xi_{0}}\xi(\t)d\t,
\end{eqnarray*}
\begin{eqnarray*}
D_{\xi_{0}}\xi(t)&=& D_{\xi_{0}}\xi(0) - \int_{0}^{t}[D_{x}^{2}a(\t) D_{\xi_{0}}x(\t)+D_{x}D_{\xi}a(\t) D_{\xi_{0}}\xi(\t)]dw(\t) \\
   &-&\frac{1}{2}\int_{0}^{t}[D_{\xi}D_{x}^{3}a(\t)D_{\xi}a(\t)+ D_{x}^{2}a(\t)D_{x}D_{\xi}a(\t) \\
  &-& D_{\xi}D_{x}^{2}a(\t)D_{x}a(\t)-D_{\xi}D_{x}a(\t)D_{x}^{2}a(\t)] D_{\xi_{0}}x(\t)d\t\\
    &-&\frac{1}{2} \int_{0}^{t}[D_{\xi}D_{x}^{2}a(\t)D_{\xi}a(\t) D_{x}^{2}a(\t)D_{\xi}^{2}a(\t) \\
   &-& D_{\xi}^{2}D_{x}a(\t)D_{x}a(\t)-(D_{\xi}D_{x}a(\t))^{2}] D_{\xi_{0}}\xi(\t)d\t.
\end{eqnarray*}
Let $n\geq 1$. It is well known that $\E |D_{\xi_{0}}x(0)|^{n} \leq K_{1}, \E |D_{\xi_{0}}\xi(0)|^{n} \leq K_{2}$, for some constants $K_{1}, K_{2}$, see Kunita\cite{K1}.
\\
We set: $\phi(t):=E\sup_{\t\leq t}|D_{\xi}x(\t)|^{n}$ and $\psi(t):=E\sup_{\t\leq t}|D_{\xi}\xi(\t)|^{n}$. Using the estimates of $D^{\alpha}_{\xi}D^{\beta}_{x}a(t)$, martingale and H\"older inequalities (we could also use the It\"o formula for the expressions of $|D_{\xi}x(t)|^{n}, |D_{\xi}\xi(t)|^{n}$), it follows from the above equations that:
\begin{equation}
\label{eq521}
\phi(t)\leq K_{n} + C_{n}\int_{0}^{t}(\phi(\t)+(1+|\xi|)^{-n}\psi(\t))d\t, 
\end{equation}
\begin{equation}
\label{eq522}
\psi(t)\leq K_{n}'+ C_{n}\int_{0}^{t}(\psi(\t)+\phi(\t))d\t,
\end{equation}
By the Gronwall lemma used in (\ref{eq522}) we get: $\psi(t)\leq K_{n}' T \phi(t) e^{C_{n}T} =:C_{n}' \phi(t)$. We substitute this in (\ref{eq521}) and 
use the Gronwall Lemma again for this inequality, we get: $\phi(t)\leq K_{n} \exp [C_{n}'((1+|\xi|)^{-n}+1)T]$, Hence $\phi(t)$ and at the same time $\psi(t)$ 
are bounded.
\\
\noindent
Next, let $\alpha\geq 2$ and suppose that for each $\gamma\leq \alpha$ we have
$E|D_{\xi}^{\gamma}x(t)|\leq C, E|D_{\xi}^{\gamma}\xi(t)| \leq C $. Then we can write:
\begin{eqnarray*}
D_{\xi}^{\alpha+1}x(t)&=& D_{\xi_{0}}^{\alpha+1}x(0) + \int_{0}^{t}[D_{x}D_{\xi}a(\t) D_{\xi}^{\alpha+1}
x(\t)+D_{\xi}^{2}D(\t)D_{\xi}^{\alpha+1}\xi(\t)]dw(\t)\\
   &+&\frac{1}{2}\int_{0}^{t}[ D_{\xi}D_{x}^{2}a(\t) D_{\xi}a(\t)+(D_{\xi}
    D_{x}^{2}a(\t))^  {2}\\
   &-& D^{2}_{\xi}D_{x}a(\t)D_{x}a(\t)-D^{2}_{\xi}D_{x}^{2}
      a(\t)]D_{x}^{\alpha+1} x(\t) d\t\\
   &+&\frac{1}{2} \int_{0}^{t}[D_{\xi}^{2}D_{x}a(\t) D_{\xi}a(\t)+D_{\xi}D_{x}
      a(\t) D_{\xi}^{2}a(\t)\\
   &-&D^{3}_{\xi}a(\t)D_{x}a(\t)-D^{2}_{\xi}D_{x}D_{\xi}a(\t)]
   D_{\xi}^{\alpha+1}\xi(\t)d\t\\
   &+& \int_{0}^{t}k_{1}(\t)dw(\t)+\int_{0}^{t}k_{2}(\t)d\t,
\end{eqnarray*}
and:
\begin{eqnarray*}
D_{\xi}^{\alpha+1}\xi(t)&=& D_{\xi_{0}}^{\alpha+1}\xi(0)-\int_{0}^{t}[D_{x}^{2}a(\t)D_{\xi}^{\alpha+1}x(\t)
+D_{x}D_{\xi}a(\t)D_{\xi}^{\alpha+1}\xi(\t)]dw(\t) \\
   &-&\frac{1}{2}\int_{0}^{t}[D_{\xi}D_{x}^{3}a(\t)D_{\xi}a(\t)+
  D_{x}^{2}a(\t)D_{x}D_{\xi}a(\t) \\
  &-& D_{\xi}D_{x}^{2}a(\t)D_{x}a(\t)-D_{\xi}D_{x}a(\t)D_{x}^{2}a(\t)]
   D_{\xi}^{\alpha+1}x(\t)d\t\\
    &-&\frac{1}{2} \int_{0}^{t}[D_{\xi}^{2}D_{x}^{2}a(\t)D_{\xi}a(\t)+
D_{x}^{2}a(\t)D_{\xi}^{2}a(\t) \\
   &-& D_{\xi}^{2}D_{x}a(\t)D_{x}a(\t)-(D_{\xi}D_{x}a(\t))^{2}]
D_{\xi}^{\alpha+1}\xi(\t)d\t\\
   &+& \int_{0}^{t}k_{1}(\t)dw(\t)+\int_{0}^{t}k_{2}(\t)d\t,
\end{eqnarray*}
where $h_{i},k_{i}$ satisfies $E |k_{i}(\t)|^{n}\leq C_{n}$, $E |h_{i}(\t)|^{n} \leq C_{n}'$. Using
the same arguments as for $\alpha=1$ we get the boundedness of $E|D_{\xi}^{\alpha+1}x(t)|^{n}$ and $E|D_{\xi}^{\alpha}\xi(t)|^{n}$.

\noindent
$\bullet$ Now, we shall show how to get the estimates for $E|D_{\xi}^{\alpha}
\bar{x}(t)|, E|D_{\xi}^{\alpha}\bar{\xi}(t)|$. Let us rewrite Eq. $({\cal C})$
 as $dX(t)=f(t, X(t))\circ dw(t)$ and consider the associated flow of 
diffeomorphisms $\phi_{s,t}(Y)$ which satisfies:
\[ \phi_{s,t}(Y)=Y-\int_{s}^{t}f(\t,\phi_{s,\t}^{-1}(Y))\circ \hat{d}w(\t). \]
Then, by Kunita [\cite{K2}, Theorem 7.3], the inverse map $\phi_{s,t}^{-1}$ 
satisfies the backward Stratonovich equation
\[ \phi_{s,t}^{-1}(Y)=Y+\int_{s}^{t}f(r,\phi_{r,t}(Y))\circ dw(r). \]
Hence, by the same arguments as above, one can prove that $E|D_{\xi}^{\alpha}
\bar{x}_{s,t}|, E|D_{\xi}^{\alpha}\bar{\xi}_{s,t}|$ satisfy the estimates
of the proposition, where $(\bar{x}_{s,t},\bar{\xi}_{s,t}):=\phi_{s,t}^{-1}$
for each $0\leq s\leq t\leq T$, which completes the proof because $(\bar{x}(t),
\bar{\xi}(t))=\phi_{0,t}^{-1}$. 
\\
\noindent
$\bullet$ If $a(t,x,\xi)$ is homogeneous of degree $1$ in $\xi$, then the fact that the solution $\phi_{t}(x_{0},\xi_{0})$ to ${\cal C}$ is homogeneous 
of degree $1$ in $\xi_{0}$ can be seen as follows: $\xi(t)$ verifies the equation:
\[ d\xi(t)=- \xi(t)\frac{\partial a}{\partial x}(t,x(t),1)\circ dw(t) \]
whose solution is of the form $\xi(t)=\xi(0)\exp[\int_{0}^{t} -\partial_{x} a (\t,x(\t),1)\circ dw(\t)].$ $\Box$
\\
\\
\noindent
{\it Proof of Corollary \ref{coroRS1}.}\\
\noindent
First, we note that the flow $\phi_{t}(x,\xi))$ is homogeneous of degree $0$ and $1$ in $x$ and $\xi$ respectively. 
Next we shall prove that given a compact set $K\in \R^{d}$ there exist constants $C(\alpha, \beta, K, w)$ such that: $\forall x \in K, \xi \in \R^{d}: |D_{\xi}^{\alpha}D_{x}^{\beta}q(x,\xi, w)|\leq C(\alpha, \beta, K, w) (1+|\xi|)^{m-|\alpha|} \; \; a.e. $
We shall apply Lemma \ref{RSlemma2} to $f(\xi,x)= D_{\xi}^{\alpha}D_{x}^{\beta}q(x,\xi, w)$. As $q(x,\xi, w)$ is homogeneous of degree $m$ in $\xi$ we have that $f(\xi,x)$ is homogeneous of degree $m-\alpha$ in $\xi$. Now we have to show that there is a constant $\gamma >0$ such that $\E|\partial_{\xi_{j}}f(\xi,x)|^{\gamma}, j=1, \cdots, d$ are
bounded. We take $\gamma = d+1$.
Given the expression of $f$ we see that the derivatives $\partial_{\xi_{j}}f(\xi,x)$ contain terms of the form:
 \[  S_{\alpha', \beta'} \times T_{\alpha", \beta"}:= 
           \frac{\partial_{\xi}^{\alpha'}\partial_{x}^{\beta'}p}{\partial \xi^{\alpha'}\partial x^{\beta'}} (\phi_{t}(x,\xi)) \times 
					\frac{\partial_{\xi}^{\alpha"}\partial_{x}^{\beta"}\phi_{t}}{\partial \xi^{\alpha"}\partial x^{\beta"}} (x,\xi) \]
To estimate $\E|\partial_{\xi_{j}}f(\xi,x)|^{\gamma}$, by using Schwarz's inequality we will have to show that the quantities: $\E|S_{\alpha', \beta'}|^{n}$ and 
$\E|T_{\alpha", \beta"}|^{n}$ are bounded with $n=2\gamma$ when $|\xi|=1$ and $x\in K$.
\\
For the terms $S$, we note that since $p\in S^{m}_{h}$, we have $|\partial_{\xi}^{\alpha'}\partial_{x}^{\beta'}p (x',\xi')| \leq C(\alpha', \beta') (1+|\xi|)^{m-|\alpha'}$ and therefore:
$ S_{\alpha', \beta'} \leq C(\alpha', \beta') (1+|\xi_{t} (x,\xi)|)^{m-|\alpha'|} $ where $(x_{t},\xi_{t})=\phi_{t}(x,\xi))$. According to Proposition \ref{propRS3} we have:
$ \E|\xi_{t} (x,\xi)|)^{n'} \leq C |\xi|^{n'}$, hence $ \E| S_{\alpha', \beta'}|^{n'} \leq C_{2}(\alpha', \beta')(1+ |\xi|)^{n'(m-|\alpha'|)} $ for some constant $C_{2}$, which implies that:
$$  \E| S_{\alpha', \beta'}|^{n'} \leq C_{3}:= C_{2}(\alpha', \beta')2^{n'(m-|\alpha'|)}  \; \; \forall \xi: |\xi|=1  $$
As for the terms $T$, we apply Proposition \ref{propRS3} and there is a constant $C_{n}'$ such that: 
$$\E |\partial_{\xi}^{\alpha"}\partial_{x}^{\beta"}\phi_{t} (x,\xi)|^{n'} \leq C_{n'}'(\alpha", \beta") (1+ |\xi|)^{n'(m-|\alpha"|)},  $$
which implies that $ \E |T_{\alpha", \beta"}|^{n'} \leq C_{n'}'(\alpha", \beta") 2^{n'(m-|\alpha"|)} $ for $ |\xi|=1  $. We can therefore apply Lemma \ref{RSlemma2} to $f$ with $\gamma=d+1$ to deduce that $q_{t} \in RS^{m}_{h}$. $\Box$


\noindent 
{\it E-mail}: adnan.aboulalaa@polytechnique.org
\end{document}